\documentclass{amsart}

\usepackage{amsmath,amscd,amsthm,amsfonts,latexsym,amssymb,theoremref,enumitem}

\usepackage{float}

\usepackage[table]{xcolor}

\usepackage{makecell}

\newcounter{countrem}
\newcounter{countt}
\newcounter{countd}
\newcounter{countp}
\newcounter{countc}
\newcounter{countl}
\newtheorem{Theorem}[countt]{Theorem}
\newtheorem{Definition}[countd]{Definition}
\newtheorem{Proposition}[countp]{Proposition}
\newtheorem{Corollary}[countc]{Corollary}
\newtheorem{Lemma}[countl]{Lemma}
\newtheorem{Remark}[countrem]{Remark}

\providecommand {\C} {{\mathbb C}}

\providecommand {\D} {{\mathcal D}}
\providecommand{\wed}{\wedge}
\providecommand{\al}{\alpha}
\providecommand{\DS}{\displaystyle}
\providecommand{\Ja}{J_a}
\providecommand{\Jd}{J_d}
\providecommand{\Jst}{\check{J}}
\providecommand{\Jast}{\check{J}_a}
\providecommand{\Jdst}{\check{J}_d}
\providecommand{\Hc}{\mathcal H}
\providecommand{\om}{\omega}
\providecommand{\ep}{\epsilon}
\providecommand{\La}{\Lambda}
\providecommand{\lamub}{\lambda_\mub}
\providecommand{\ladeb}{\lambda_\delb}
\providecommand{\lade}{\lambda_\del}
\providecommand{\lamu}{\lambda_\mu}
\providecommand{\rhmub}{\rho_\mub}
\providecommand{\rhdeb}{\rho_\delb}
\providecommand{\rhde}{\rho_\del}

\providecommand{\tadeb}{\tau_\delb}
\providecommand{\tade}{\tau_\del}
\providecommand{\tamub}{\tau_\mub}
\providecommand{\tamu}{\tau_\mu}
\providecommand{\grcel}{\cellcolor{gray!25}}
\providecommand\del{\partial}
\providecommand\delb{{\bar{\partial}}}
\providecommand\mub{{\bar{\mu}}}
\providecommand\clM{{{\rm Cl}(M)}}
\providecommand\ClM{{{\C}{\rm l}(M)}}

\providecommand{\iprod}{\mathbin{\lrcorner}}

\makeatletter
\@namedef{subjclassname@2020}{%
  \textup{2020} Mathematics Subject Classification}
\makeatother

\title[K\"ahler identities for almost complex manifolds]
{K\"ahler identities for almost complex manifolds}
\author{Luis Fern\' andez and Samuel Hosmer}

\begin{document}

\begin{abstract}
  We obtain a generalization, for a general compact almost complex manifold,
  of the well-known K\"ahler (or Hodge) identities for K\"ahler manifolds involving the commutators of the exterior differential and the Lefschetz operator and its adjoint. The main idea is to study the problem on the Clifford bundle via the Dirac operator,
  and then translate the results to the exterior bundle.
\end{abstract}
\sloppy
\maketitle

\section{Introduction}

On a K\"ahler manifold, the K\"ahler (or Hodge) identities give relationships between the commutator of the exterior differential $d$ and the Lefschetz operator on forms given by $L(\al) = \omega\wed\al$, where $\omega$ is the fundamental form. Namely, we have
$$[d^*, L] = -d^c \ \ \ \qquad\qquad \ \ \ [d,L] = 0,$$
where the bracket denotes the commutator, $d^* := -*\circ d \circ *$ is the formal adjoint of $d$, and $d^c := \Ja^{-1}\circ d\circ \Ja = i(\delb-\del)$, with $\Ja$ being the extension of the complex structure $J$ as an algebra operator on the exterior algebra. These identities seem to have appeared first in \cite{hodgepaper, hodgebook}, and in modern notation in \cite{weilbook}. They are a fundamental tool to prove several important results for K\"ahler manifolds, such as the Lefschetz decomposition of complex DeRham cohomology or the fact that, on a K\"ahler manifold, the notions of $d$-harmonic, $\del$-harmonic, and $\delb$-harmonic forms coincide.

There are generalizations of the K\"ahler-Hodge identities: In \cite{demaillytau}, Demailly defines new operators
\begin{equation*}
  \lambda(\alpha) = d\omega\wedge\alpha  \ \ \mbox{and} \ \ \tau = [\Lambda, \lambda],
\end{equation*}
and shows that, on a {\em complex} manifold (i.e. integrable), the K\"{a}hler identities extend to
\begin{equation}
  [d^*,L] = -d^c -\tau^c  \ \ \ \qquad\qquad \ \ \ [d,L] = \lambda.\label{dema}
\end{equation}
For {\em almost K\"ahler} manifolds (i.e. non-integrable almost complex manifolds with a closed fundamental form), similar identities are proved in \cite{bartolo} and \cite{ciriciwilson}. The proofs and generalizations we know of these identities rely on a local analytic study of the problem.

In this paper we prove the identities (\ref{dema}) for general {\em almost complex} manifolds (i.e. just a manifold with an almost complex structure $J$).
Our approach does not rely on a local analytic study of the problem. Instead, we formulate the problem in the Clifford bundle of the manifold instead of the exterior bundle, following the ideas used in \cite{Mich} for K\"ahler manifolds, and then we translate the results to the exterior algebra. This offers new possibilites of study due to the somewhat richer algebraic structure of the Clifford algebra, and we think that there still a lot to explore using this formulation.

We make an effort to distinguish between the Clifford bundle and the exterior bundle, even though they are isomorphic via the metric. The aim is to emphasize the richer structure of the former and to separate the idea of a ``Clifford cohomology'' from the usual DeRham cohomology. While this separation requires additional notation, we believe that it makes the exposition more clear.

The paper is organized as follows: in Section 2 we make the necessary definitions, present the basic preliminary results, and fix the notation that will be used throughout the paper. In Section 3 we formulate the problem of finding an equivalent of the K\"{a}hler identities in the Clifford bundle, and we prove the Clifford version of these identities, which is surprisingly simple. We also prove other interesting commutation identities. Finally, Section 4 is devoted to translate these identites to the exterior algebra. This translation becomes rather technical and often requires computations using a local frame, but it is nevertheless a linear algebraic procedure.

\section{Preliminaries and Notation}\label{notation}

Let $M$ be a $2n$ dimensional manifold with an almost complex
structure, that is, a smooth bundle map $J:TM\to TM$ satisfying $J^2 =
-I$. We will use the same notation for the complex-linear extension 
$J: T_\C M\to T_\C M$, and will denote by $\Jst$ the dual map $\Jst: T^*_\C M\to  T^*_\C M.$ In general, the dual of an operator $P$ will be denoted by $\check{P}$, and its adjoint by $P^*$.

Let $g$ be a metric compatible with $J$, that is, $\langle X,Y\rangle = \langle JX, JY\rangle$. Then we can construct the Clifford bundle $\clM$ over $TM$, and its complexification $\ClM$. With a slight abuse of notation, the space of smooth sections over these bundles will also be denoted by $\clM$ and $\ClM$, respectively. Likewise, $A(M)$ (or $A_\C(M)$) will denote both the exterior bundle (or complexified exterior bundle) and the space of smooth forms (or complex smooth forms).

The sesquilinear extension of the metric $g$ to the complexification
will be also be denoted by $g$. As is well known,
$\clM$ ($ = {\rm Cl}(TM)$) and $\Lambda(TM)$ are canonically isomorphic (see for
example \cite{LawMich}), 
 and by extension so are $\ClM$ and  $\Lambda_{\C}(TM)$.
In turn, via the metric, $\Lambda_\C(TM)$ is
isomorphic to $A_\C(M) \equiv\Lambda_{\C}(T^*M)$ via the map
$$X\longrightarrow (Y\to \langle Y, \bar{X}\rangle),$$
where the bar denotes conjugation. The composition
of the two isomorphisms will be denoted by
$$\flat : \ClM\to A_\C(M), \ \ \ \mbox{with inverse} \ \ \ \sharp:=\flat^{-1}:
A_\C(M)\to \ClM.$$
We will often denote $\flat(X)$ by $X^\flat$ and $\sharp(\al)$ by
$\al^\sharp$. Note that, if $\{e_A\}_{A=1}^{2n}$ is an orthonormal basis with dual basis $\{\theta^A\}_{A=1}^{2n}$, then $\flat(e_{A_1}\cdot
  \ldots\cdot e_{A_k}) = \theta^{A_1}\wedge \cdots\wedge \theta^{A_k} $. We will often use the symbol ``$\cong$'' to denote 
correspondence via these ``musical'' isomorphisms, so for example we could have writen $e_{A_1}\cdot
  \ldots\cdot e_{A_k} \cong \theta^{A_1}\wedge \cdots\wedge \theta^{A_k} $ instead.

The maps $J:TM\to TM$ and $\Jst:T^*M\to T^*M$ can be extended to $\ClM$
  and $A_\C(M)$, respectively, in two different ways:
  \begin{itemize}[leftmargin=20pt]
  \item As a derivation $\Jd$:
    $\Jd(X\cdot Y) = JX\cdot Y + X\cdot JY$ for $X, Y\in T_\C M$
    (and similarly for $\Jdst$).
  \item As an algebra map $\Ja$: $\Ja(X\cdot Y) = JX\cdot JY$ for
    $X, Y\in T_\C M$ (and similarly for $\Jast$).
  \end{itemize}
Notice that $\Ja$ and $\Jd$ commute, and so do $\Jast$ and $\Jdst$.
  
The algebra $A(M)$ is $\mathbb{Z}$-graded and $\clM$ is
$\mathbb{Z}_2$-graded. The eigenspaces of $\Jast$ in each degree induce a $\mathbb{Z}\times \mathbb{Z}$ bigrading of $A_{\C}(M)$; the algebra $\ClM$ is also 
$\mathbb{Z}\times \mathbb{Z}$-bigraded (see \cite{Mich}).    Given a complex form $\psi$, its $(p,q)$ part in the bigrading $A_\C(M) = \oplus A^{p,q}_\C(M)$ will be denoted by $\psi^{p,q}$. Also, if $\psi\in A^3_\C(M)$,  we will use the same notation as in \cite{gauduchon}:
   $$\psi^+ := \psi^{1,2} + \psi^{2,1} \ \ \ \mbox{and} \ \ \ \psi^- := \psi^{0,3} + \psi^{3,0}.$$
An operator $P$ is said to have {\em degree
$p$}, or {\em bidegree
$(p,q)$}, if it takes elements of degree $r$ to elements of degree
$r+p$, or elements of degree $(r,s)$ to elements of degree
$(r+p,s+q)$, respectively.

Given operators $P$, $Q$ in $\ClM$ or $A_{\C}(M)$ of total
  degree $p$ and $q$, respectively, $[P , Q ]$ will denote the {\em
    supercommutator} of $P$ and $Q$, that is,
$$[P,Q] = P Q - (-1)^{pq} QP.$$ 
Notice that $[P,Q] = -(-1)^{pq} [Q,P]$, and that if $P$ and $Q$ are derivations, then $[P,Q]$ is also a derivation. The supercommutator satisfies the Jacobi identity
$$[P,[Q,R]] = [[P,Q],R] + (-1)^{pq} [Q,[P,R]].$$

    As usual, given an operator $P$ in $\ClM$ or $A_{\C}(M)$, the superindex $c$ will denote conjugation with
  $\Ja$ (or $\Jast$), that is
$$P^c = \Ja^{-1}\circ P\circ \Ja \ \mbox{(or $P^c = \Jast^{-1}\circ
  P\circ \Jast).$}
$$
Notice that if $P$ is an operator on $A_C(M)$ of bidegree $r,s$, then if $\phi\in A^{p,q}$,
$$P^c(\phi) = \Jast^{-1}(P(\Jast\phi)) = i^{p-q}\Jast^{-1}(P(\phi)) = i^{p-q} i^{q+s-p-r}P(\phi) = i^{s-r} P(\phi),$$
and therefore
\begin{equation}
  P^c = i^{s-r} P \ \ \ \mbox{if $P$ has bidegree $(r,s)$}\label{Pcdeg}
\end{equation}
All these objects satisfy the following elementary properties that will be used throughout.

\begin{Lemma}\thlabel{elemprop}
(Elementary properties.) Let $\psi\in A^k_\C (M)$.
\begin{enumerate}[label=\rm{(\alph*)},leftmargin=25pt]
\item $\flat\circ J\circ\sharp = -\Jst$.
  \label{2a}
\item $(J_a)^\vee = (\Jst)_a$ (that is, extending as an algebra map
commutes with duality).\label{2b}
\item $(J_d)^\vee = (\Jst)_d$ (that is, extending as a derivation
  commutes with duality).\label{2bb}
\item $\flat\circ\Jd\circ\sharp = -\Jdst$.
\label{2d}
\item $\flat\circ\Ja\circ\sharp (\psi) = (-1)^k \Jast
  \psi$. \label{2c}
\item $(\Jst_a)^* \psi = (-1)^k \Jst_a\psi$ and $\Jdst^* = -\Jdst$.
  (and the same identities for $\Ja$ and $\Jd$).
  \label{2f}
\item $(\Jst_a)^{-1} \psi = (-1)^k \Jst_a\psi$.
  (and the same identity for $\Ja$).
  \label{2g}
\item If $P$ is an operator of degree $p$ on $\Gamma\ClM$ or
    $A_\C(M)$, then $(P^c)^* = (P^*)^c$ and
    $(P^c)^c = (-1)^p P$.\label{2h}
\item If $P:\Gamma\ClM\to \Gamma\ClM$ is an operator of degree $p$
    then \
    $\flat\circ P^c\circ\sharp = (-1)^p (\flat\circ P\circ
    \sharp)^c$. \label{Pequiv}
\item If $P:\Gamma A_\C(M)\to \Gamma A_\C(M)$ is an operator of {\em
      bidegree} $(p,q)$ then $P^c = i^{q-p} P$. \label{Pc}
\end{enumerate}
\end{Lemma}
\begin{proof} They are all elementary consequences of the definitions. For example, to prove \ref{2a}, let $X\in T_\C(M)$ and $\beta\in T^*_\C(M)$. Since $\flat$
  and $\sharp$ are isometries,
$$\langle X,J \beta^\sharp\rangle= \langle JX,-\beta^\sharp\rangle =
\langle (JX)^\flat,-\beta\rangle = -\beta(JX) = -\Jst\beta(X) =
\langle X, (-\Jst\beta)^\sharp\rangle,$$ so
$J \beta^\sharp = (-\Jst\beta)^\sharp$, or
$\flat\circ J\circ\sharp ( \beta) = -\Jst\beta$.

\end{proof}

The fundamental form $\omega\in A^2(M)$ is defined by
$\omega(X,Y) := \langle JX, Y\rangle$. In terms of a local adapted
orthonormal coframe, i.e. an orthonormal coframe
$\{\theta^A\}_{A=1}^{2n}$ such that $\Jst\theta^i = -\theta^{i+n}$, $1\le
i\le n$, 
$\omega$ can be written as
$$\omega = \sum_{i=1}^{n} \theta^i\wed\theta^{i+n}$$
(in general,
   uppercase sub- and superindices will run from 1 to $2n$, whereas lowercase
   sub- and superindices will run from $1$ to $n$).\\
   The {\em Lee form} $\theta$ is defined by (see for example \cite{gauduchon}).
   \begin{equation}
     \theta:= \omega\iprod d\omega^+ = -\Jst d^*\omega^+ \label{leeform2}
   \end{equation}
For forms $\alpha, \beta$, we write $\alpha \iprod\beta$ to denote the adjoint of $\beta\to \alpha\wedge\beta$, or equivalently, contraction of $\beta$ with $\alpha^\sharp$.

   The same symbol $\omega$ will be used to denote $\omega^\sharp\in\clM$, which can be expressed, in terms of a frame, by
   $$\omega = \sum_{i=1}^{n} e_i\cdot e_{i+n}.$$

Following \cite{Mich}, we define the following operators in $\ClM$:

\begin{Definition}\thlabel{firstdef}
 Let $\nabla$ denote the Levi-Civita connection on $TM$ or its extension to  $\ClM$ or $A_\C(M)$.
 \begin{itemize}[leftmargin=20pt]
 \item The Dirac operator $D:\ClM\to\ClM$ is defined by
   $$DX = \sum_{A=1}^{2n} e_A\cdot \nabla_{e_A} X,$$
   where $\{e_A\}_{A=1}^{2n}$ is an orthonormal basis (it is immediate to check
   that the definition is independent of the basis chosen).
 \item $\Hc:\ClM\to\ClM$ defined by (see \cite{Mich}, p. 1091)
$$\Hc ( X )  =\dfrac{1}{2i}\left(\om\cdot  X  +  X \cdot\om
\right).$$
\item For $X\in \ClM$, $L_X(Y):= X\cdot Y$ and $R_X:= Y\cdot X$. 
\end{itemize}
\end{Definition}

These objects have the following basic properties.

\begin{Lemma}\thlabel{elemprop2}
We have
\begin{enumerate}[label=\rm{(\alph*)},leftmargin=25pt]
\item $J_d\omega = 0$ and $J_a\omega = \omega$.\label{elp1}
\item For $ X \in \ClM$,
    $\Jd( X ) = \frac{1}{2}(\omega\cdot X - X \cdot \omega).$
    Therefore, $\Hc = i \Jd - i L_\om$.\label{wJ}
\item $\nabla_X \Jd$ and $\Ja^{-1}\circ\nabla_X \Ja$ are derivations
    of order 0, and they are anti-self adjoint. \label{JaJd}
\item $\Hc J_d = J_d \Hc$, $\Hc^c = \Hc$, and $\Jd^c =\Jd$ (i.e. $J_d$ and $J_a$ commute).\thlabel{Hc}
\item $\nabla \psi^\sharp = (\nabla\psi)^\sharp$ (that is,
    $\flat\circ\nabla\circ\sharp = \nabla$).\label{2j}
\item $\Ja\nabla\omega = -\nabla\omega$ and
    $\Jast\nabla\omega = -\nabla\omega$.\label{2n}
\end{enumerate}
\end{Lemma}

\begin{proof}
For \ref{elp1}, note that for an adapted frame, if $1\le i\le n$,  $J_d(e_i\cdot e_{i+n}) = e_{i+n}\cdot e_{i+n} - e_i\cdot e_i = 0$ and $J_a(e_i\cdot e_{i+n}) = - e_{i+n}\cdot e_{i} = e_i\cdot e_{i+n}.$

For \ref{wJ}, see \cite{Mich}, Lemma 2.5 (note though that our
definition of $J$ is $i$ times their definition of $J$).

To show \ref{JaJd}, by definition,
$(\nabla_X \Jd)(Y) = \nabla_X (\Jd Y) - \Jd \nabla_X Y = [\nabla_X, \Jd]
(Y)$. Therefore $\nabla_X\Jd$ is a derivation since it is the commutator
of two derivations. To show that $\Ja^{-1}\circ(\nabla_X \Ja)$ is a
derivation, note that $\Ja^{-1}\circ(\nabla_X \Ja) = \Ja^{-1}\circ \nabla\circ  \Ja - \nabla$, which is a sum of derivations and therefore also a derivation.
To show that $\nabla_X\Jd$ is anti-self adjoint, let
$X\in T_\C M$, $Y, Z\in \ClM$, recall that $\Jd^* = -\Jd$ (\thref{elemprop}\ref{2f}), and
compute
\begin{eqnarray*}
  \langle (\nabla_X\Jd) Y, Z\rangle 
  &=& \langle \nabla_X (\Jd Y), Z\rangle - \langle \Jd \nabla_X Y,
      Z\rangle\\
  &=& X\langle\Jd Y, Z\rangle - \langle \Jd Y, \nabla_X
      Z\rangle + \langle  \nabla_X Y, \Jd Z\rangle\\
  &=& -X\langle Y, \Jd Z\rangle + \langle  Y, \Jd\nabla_X
      Z\rangle + X\langle  Y, \Jd Z\rangle -\langle  Y, \nabla_X (\Jd Z\rangle)\\
  &=& -\langle Y,(\nabla_X \Jd) Z\rangle.
\end{eqnarray*}
A similar computation shows that $\Ja^{-1}\nabla J_a$ is anti-self adjoint.

Item \ref{Hc} follows directly from \ref{elp1}, the definition of $\Hc$, and \ref{wJ}.

Part \ref{2j} is essentially the definition of $\nabla$ in $A_\C(M)$.

To show \ref{2n}, use the definition of $\omega\in \Gamma A^2_\C (M)$
and the properties of $\nabla$. Let $X\in T_\C M$,
and let $Y,Z$ be sections of $T_\C M$. Then
\begin{eqnarray*}
  (\Jast\nabla_X\omega) (Y,Z) &=& \nabla_X\omega (JY,JZ)\\
                              &=& X(\omega(JY,JZ)) 
                                  - \omega(\nabla_X JY,JZ)
                                  -\omega(JY, \nabla_X JZ)\\
                              &=&  
                                  -X\langle Y,JZ\rangle 
                                  - \langle J\nabla_X JY,JZ\rangle
                                  +\langle Y, \nabla_X JZ\rangle\\
                              &=&  
                                  -X\langle Y,JZ\rangle 
                                  - X\langle  JY,Z\rangle
                                  +\langle JY, \nabla_X Z\rangle
                                  + X\langle Y,JZ\rangle 
                                  -\langle \nabla_X Y, JZ\rangle\\
                              &=&
                                  -(X\omega(X,Y) - \omega(Y,\nabla_X Z) -\omega(\nabla_X Y,Z))\\
                              &=&-(\nabla_X\omega)(Y,Z)
\end{eqnarray*}
The equivalent statement for $\omega\in\Gamma\ClM$ follows then from
Lemma \thref{2d} and \ref{2j}.

\end{proof}

\section{K\"ahler identities in the Clifford bundle}

The classical K\"ahler identities read
$$[d^*, L] = -d^c \ \ \ \qquad\qquad \ \ \ [d,L] = 0  \ \ \ \qquad\qquad \ \ \ [d, \Lambda] = (d^c)^*  \ \ \ \qquad\qquad \ \ \ [d^*,\Lambda] = 0,$$
where $L(\phi) := \omega\wedge\phi$ is the Lefschetz operator and $\Lambda(\phi):= \sum_{i=1}^n \omega \iprod \phi$ its adjoint.
In order to find an equivalent formulation in the Clifford bundle, we need the following conversion between operators.

\begin{Proposition}\thlabel{equiv1}
  There are the following equivalences between operators in $\ClM$ and
  operators in $A_\C(M)$.
  \begin{enumerate}[label=\rm{(\alph*)},leftmargin=25pt]
  \item $D\cong d+d^*$.\label{Dequiv}
  \item $\Hc\cong i(\Lambda-L)$.\label{Hcequiv}
  \item $[D, \Hc]\cong i[d+d^*,\Lambda - L]$.\label{commDHequiv}
  \item $D^c\cong - (d^c+ {d^*}^c)$.\label{Dcequiv}
  \end{enumerate}
\end{Proposition}
\begin{proof}
  For \ref{Dequiv} see for example \cite{LawMich}, Theorem 5.12. For
  \ref{Hcequiv} see \cite{Mich}, page 1118.  Part \ref{commDHequiv}
  follows immediately from \ref{Dequiv} and \ref{Hcequiv}. Part
  \ref{Dcequiv} follows immediately from \ref{Dequiv} and
  \thref{elemprop} \ref{Pequiv}, since $D$ has degree 1.

\end{proof}

Thus, as a guide, to find an equivalent formulation of the K\"{a}hler identities it seems sensible to find an expression for $[D, \Hc]$ which contains a term involving $D^c$. Miraculously, one obtains a nice expression once a suitable operator is defined. More precisely, we have

\begin{Lemma} \thlabel{DJcomm} Let $\{e_A\}$ an adapted frame as before. Then
  $$[D, \Hc] = -iD^c -i L_{D\omega} + i \sum_{A=1}^{2n} e_A\cdot (\nabla_{e_a} \Jd + \Ja^{-1}\circ \nabla_{Je_A} \Ja).$$
\end{Lemma}

\begin{proof}
  Since $\Hc = i\Jd - i L_\omega$, we first compute the commutator with each of the two terms and compare the result with $D^c$.
\begin{eqnarray*} [D, \Jd] ( X ) &=& \sum_{A=1}^{2n}\bigl( e_A\cdot
    \nabla_{e_A} (\Jd X ) -
                                       \Jd ( e_A\cdot \nabla_{e_A}  X ) \bigr)\\
                                   &=& \sum_{A=1}^{2n}\bigl( e_A\cdot
                                       (\nabla_{e_A} \Jd) X +e_A\cdot
                                       \Jd\nabla_{e_A} X -J e_A\cdot
                                       \nabla_{e_A} X
                                       -  e_A\cdot \Jd\nabla_{e_A}  X \bigr)\\
                                   &=& \sum_{A=1}^{2n} e_A\cdot
                                       (\nabla_{e_A} \Jd) X  -\sum_{A=1}^{2n} J e_A\cdot
                                       \nabla_{e_A} X.
\end{eqnarray*}
Also,
\begin{eqnarray*}
[D,L_\om](X) &=&  \sum_{A=1}^{2n}(e_A\cdot \nabla_{e_A}(\om\cdot X) - \om\cdot e_A\cdot \nabla_{e_A} X)\\
             &=&
                 \sum_{A=1}^{2n}(e_A\cdot \nabla_{e_A}\om\cdot X + e_A\cdot \om\cdot \nabla_{e_A}X - \om\cdot e_A\cdot \nabla_{e_A} X)\\
             &=& L_{\D\om} (X)  - 2\sum_{A=1}^{2n} J e_A\cdot \nabla_{e_A}X \ \ \mbox{using Lemma \thref{wJ}.}               
\end{eqnarray*}
Thus, since  $\Hc = i\Jd - i L_\omega$,
$$[D,\Hc] + i L_{\D\om} (X) = i \sum_{A=1}^{2n} e_A\cdot
(\nabla_{e_A} \Jd) X + i \sum_{A=1}^{2n} J e_A\cdot \nabla_{e_A}X .$$
On the other hand,
\begin{eqnarray*}
  D^c(X) &=& \Ja^{-1}\left(\sum_{A=1}^{2n} e_A\cdot \nabla_{e_A} \Ja X\right)\\
         &=& \sum_{A=1}^{2n} \left(-Je_A\cdot \Ja^{-1}((\nabla_{e_A} \Ja) X + \Ja\nabla_{e_A} X)\right)\\
         &=&  \sum_{A=1}^{2n} e_A\cdot \Ja^{-1}(\nabla_{Je_A} \Ja) X - \sum_{A=1}^{2n} Je_A\cdot\nabla_{e_A} X,\\
\end{eqnarray*}
where the first term in the last equality comes from the fact that if $1\le A\le n$, then $Je_A = e_{A+n}$, and if $n+1\le A\le 2n$, then $Je_A = e_{A-n}$. 
Therefore,
$$[D, \Hc] + i L_{\D\om}   + iD^c  = i \sum_{A=1}^{2n} e_A\cdot (\nabla_{e_a} \Jd + \Ja^{-1}\circ \nabla_{Je_A} \Ja),$$
which proves the claim.

\end{proof}

This motivates the definition of the following operators.

\begin{Definition} Let
  \begin{enumerate}[label=\rm{(\alph*)},leftmargin=25pt]
  \item For $X\in T_\C M$, define
    $\sigma_X, 
    : \ClM\to \ClM$ by
    $$\sigma_X = \nabla_X \Jd + \Ja^{-1}\circ\nabla_{JX} \Ja
    .$$
   \item For $Y\in \ClM$, let 
     $$D_\sigma = \sum_{A=1}^{2n} e_A\cdot \sigma_{e_A}( Y)
     .$$
\end{enumerate}
\end{Definition}

We will see below that $\sigma$ and $D_\sigma$ have remarkable properties. But first we have

\begin{Theorem} \thlabel{Dcomm} {\rm (Generalized K\"ahler identities
    in the Clifford bundle.)}
$$[D, \Hc] = -i D^c +i D_\sigma - i L_{D\om}.$$
\end{Theorem}
\begin{proof}
  Immediate from the previous lemma and definition.

\end{proof}

In the next section we will see that the conversion of this identity to the exterior bundle gives a generalization of the classical identities. Next we study some of the properties of the operators $\sigma_X$ and $D_\sigma$ defined above.

\begin{Lemma}\thlabel{sigmaprops}
  The operator $\sigma$ has the
  following properties.
  \begin{enumerate}[label=\rm{(\alph*)},leftmargin=25pt]
  \item $\sigma_X$ is a derivation of degree 0.\thlabel{sigmader}
  \item $\sigma_X$ is anti-self adjoint.\thlabel{sigmasa}
\item If $Y\in T_\C M$, then \
  $\sigma_{JX} (Y) = J\sigma_X (Y)= -\sigma_X(JY).$\label{sigmaJ}
\item $(\sigma_X)^c = -\sigma_X$.\label{sigmac}
\item $[\sigma_X, \Jd] = -2 \sigma_{JX}$.\label{sigmabra}
\end{enumerate}
\end{Lemma}
\begin{proof}
  (a) and (b) follow immediately from \thref{elemprop} \ref{JaJd}.

For \ref{sigmaJ}, let $X,Y\in T_\C M$ and note that, on $T_\C M$,
  $\Ja \equiv \Jd \equiv J$. Then
$$\sigma_{JX} (Y) = (\nabla_{JX} J) (Y) -J^{-1}(\nabla_{X} J)(Y) = 
J ( (\nabla_{X} \Jd)(Y) + \Ja^{-1} (\nabla_{JX} \Ja) (Y)) = -J\sigma_X
Y,$$ and
\begin{eqnarray*}
  \sigma_{X} (JY) &=& (\nabla_{X} \Jd) (JY) +\Ja^{-1}(\nabla_{JX}
                      \Ja)(JY)\\
                  &=& -\nabla_{X} Y - J\nabla_X J Y  
                      -J^{-1}\nabla_{JX} Y - \nabla_{JX} J Y \\
                  &=& -J (\nabla_X J Y -J\nabla_{X} Y
                      + J^{-1}\nabla_{JX} J Y - \nabla_{JX} Y)\\
                  &=& -J \left((\nabla_X \Jd) Y
                      + \Ja^{-1}(\nabla_{JX} \Ja) Y\right)\\
                  &=& -J\sigma_X Y.
\end{eqnarray*}
To prove \ref{sigmac}, given $Y\in T_\C M$, use \ref{sigmaJ} to obtain
  $(\sigma_X)^c(Y) = -J ( \sigma_X(JY)) = -J(-J\sigma_X Y) = -\sigma_X
  Y$. Thus, $(\sigma_X)^c = -\sigma_X$ on vectors and scalars (note
  that $\sigma \equiv 0$ on scalars), and since both sides are
  derivations of degree 0, they must agree in all of $\ClM$.
  
For \ref{sigmabra}, if $X,Y\in T_\C M$,
  $[\sigma_X, \Jd](Y) = \sigma_X(\Jd(Y)) - \Jd(\sigma_X(Y)) = -2
  \sigma_{JX}Y$ by part \ref{sigmaJ}. Thus, $[\sigma_X, \Jd]$ and
  $-2 \sigma_{JX}$ are derivations that agree on scalars and vectors,
  so they must be equal.

\end{proof}

To prove further properties of $\sigma$ we need the following technical result
that appears in
\cite{gauduchon}. It will also be used in the next section when we
translate the operator $\sigma_X$ into the exterior algebra.  
\begin{Lemma}\thlabel{psi+l} For 3-forms in an almost complex manifold
  $M$ we have:
  \begin{itemize}[leftmargin=18pt]\setlength\itemsep{6pt}
  \item If $\psi\in A^{2,1}_\C(M)\oplus A^{1,2}_\C(M)$ and
    $X,Y,Z\in T_\C M$. Then\vspace{3pt}
    \begin{enumerate}[label=\rm{(\alph*)},leftmargin=20pt]
    \item
      $\psi(X,Y,Z) = \psi(JX,JY,Z) + \psi(JX,Y,JZ) +
      \psi(X,JY,JZ)$.\vspace{5pt}\label{psi+}
    \item
      $\DS\sum_{A,B=1}^{2n} \left(\psi(e_A, Z,e_B) - \psi(e_A, JZ,
        Je_B)\right) \theta^A\wed\theta^B = \dfrac{1}{2}
      \DS\sum_{A,B=1}^{2n} \left(\psi(e_A, Z,e_B)\, + \psi(Je_A,
        Z,Je_B)\right)\, \theta^A\wed\theta^B.$\label{psi++}
    \end{enumerate}
  \item If $\xi\in A^{3,0}_\C(M)\oplus A^{0,3}_\C(M)$,
    $X,Y,Z\in T_\C M$. Then $\xi(JX,Y,Z) = \xi(X,JY,Z) = \xi(X,Y,JZ)$.
  \end{itemize}
\end{Lemma}
\begin{proof} These results can be proved by straightforward computations. For \ref{psi+}, note that if $\psi\in A^{2,1}_\C(M)\oplus A^{1,2}_\C(M)$, then $\psi^{2,1}(X,Y,Z) = \psi(X^{1,0},Y^{1,0},Z^{0,1})
  +\psi(X^{1,0},Y^{0,1},Z^{1,0})+\psi(X^{0,1},Y^{1,0},Z^{1,0})$, and similarly for $\psi^{1,2}$. Comparing this with 
$\psi(JX,JY,Z) + \psi(JX,Y,JZ) + \psi(X,JY,JZ)$ gives the result.
  
For item (b), use part (a) and compute to obtain
$$\DS\sum_{A,B=1}^{2n} \bigl(
    \psi(e_A, Z,e_B) \;\theta^A\wed\theta^B -
    \psi(Je_A, Z, Je_B) \;\theta^A\wed\theta^B\bigr)
       =   2\DS\sum_{A,B=1}^{2n} 
        \psi(e_A, JZ, Je_B) \;\theta^A\wed\theta^B.
        $$
        Then substitute the expression for $\sum_{A,B=1}^{2n} 
        \psi(e_A, JZ, Je_B) \;\theta^A\wed\theta^B$ on the right hand side of \ref{psi+}.

For (c), if $\xi\in A^{3,0}_\C(M)\oplus A^{0,3}_\C(M)$,
    $$\xi(JX,Y,Z) = \xi((JX)^{1,0},Y^{1,0},Z^{1,0}) +
                    \xi((JX)^{0,1},Y^{0,1},Z^{0,1})=
                     i\xi(X^{1,0},Y^{1,0},Z^{1,0}) -
                     i\xi(X^{0,1},Y^{0,1},Z^{0,1}),$$
and the same expression is obtained for $\xi(X,JY,Z)$ and $\xi(X,Y,JZ)$.

\end{proof}
The following properties show that $\sigma$ has a very interesting relationship with objects in the exterior algebra. In fact, it gives an idea of the $(2,1) + (1,2)$ part of the fundamental form, as follows.

\begin{Lemma}\thlabel{sigmaprops2}
  The operator $\sigma$ has the
  following properties.
  \begin{enumerate}[label=\rm{(\alph*)},leftmargin=25pt]
  \item For $Y\in T_\C M$,
$\sigma_X(Y) = \DS\sum_{B=1}^{2n} (d\omega^+(X,Y,e_B) -
d\omega^+(X,JY,Je_B))e_B.$ \thlabel{sigmaom1}
\item
  $\sigma^\flat_X (:=\flat\circ\sigma_X\circ\sharp) = -\nabla_X \Jdst
  + \Jast^{-1}\circ \nabla_{JX} \Jast$, and therefore it is a
  derivation by \thref{elemprop2} \ref{JaJd}.

  \noindent Further, if $\al$ is a
  1-form,
$\sigma^\flat_X(\al)  = 
\DS\sum_{B,C=1}^{2n} (d\omega^+(X,e_C,e_B) -
d\omega^+(X,Je_C,Je_B))\,\al(e_C)\;\theta^B.$\thlabel{sigmaflatder}
\item $\DS\sum_{A=1}^{2n} \sigma_{e_A}(e_A) = -2(\Jst \theta)^\sharp$,
  where $\theta:= \omega\iprod d\omega^+ = -\Jst d^*\omega^+$ is the
  Lee form defined in (\ref{leeform2}).\thlabel{leeform}
\end{enumerate}
\end{Lemma}
\begin{proof}
  
\noindent \ref{sigmaom1}: Using Proposition 4.2 in page 148 of \cite{KN2} we have,
adjusting the constants to our definitions, that for vectors
$X,Y,Z\in T_\C M$,
\begin{equation}\label{KN1}
 2\langle(\nabla_X J) Y,Z\rangle = d\omega(X,Y,Z) - d\omega(X,JY,JZ)
 + 4 \langle JX, N(Y,Z)\rangle,
\end{equation}
where $N$ is the Nijenhuis tensor defined by
$$N(X,Y) = \dfrac{1}{4}([JX,JY] -J[JX,Y] - J[X,JY] - [X,Y]).$$
This implies, using the fact that $N(Y, JZ) = -JN(Y,Z)$,
\begin{equation}\label{KN2}
2\langle J^{-1}(\nabla_{JX} J) Y,Z\rangle = d\omega(JX,Y,JZ) + d\omega(JX,JY,Z)
- 4 \langle JX, N(Y,Z)\rangle,
\end{equation}
and adding (\ref{KN1}) and (\ref{KN2}),
\begin{eqnarray*}
  2\langle (\nabla_X J +J^{-1}\nabla_{JX} J) Y,Z\rangle &=&
                                                            d\omega(X,Y,Z) - d\omega(X,JY,JZ) + 
                                                            d\omega(JX,Y,JZ) + d\omega(JX,JY,Z)\\
                                                        &=&
                                                            d\omega^+(X,Y,Z) - d\omega^+(X,JY,JZ) + 
                                                            d\omega^+(JX,Y,JZ) + d\omega^+(JX,JY,Z)\\ 
                                                        &\quad+& d\omega^-(X,Y,Z) - d\omega^-(X,JY,JZ) + 
                                                                 d\omega^-(JX,Y,JZ) + d\omega^-(JX,JY,Z)\\
                                                        &=&
                                                            2(d\omega^+(X,Y,Z) - d\omega^+(X,JY,JZ)) \ \ {\mbox{by \thref{psi+l}}}.
\end{eqnarray*}
Thus,
$\langle\sigma_x(Y),Z\rangle = d\omega^+(X,Y,Z) - d\omega^+(X,JY,JZ)$,
which proves \ref{sigmaom1}.

\noindent \ref{sigmaflatder}: The first identity follows from the fact that, using
  \thref{elemprop} \ref{2d}, \ref{2c} and \thref{elemprop2} \ref{2j}, if
  $\psi\in A^k_\C(M)$, then
$$\flat\circ(\nabla
\Jd)\circ\sharp (\psi) = \flat(\nabla \Jd\psi^\sharp -\Jd\nabla
\psi^\sharp) =-\nabla\Jdst\psi +\Jdst\nabla \psi =
-(\nabla\Jdst)\psi,$$ and
$$\flat\circ\Ja^{-1}\circ(\nabla \Ja)\circ\sharp (\psi) = \flat(
\Ja^{-1}(\nabla \Ja\psi^\sharp) - \nabla\psi^\sharp) = (-1)^{k+k}
\Jast^{-1}(\nabla \Jast\psi^\sharp) - \nabla\psi^\sharp =
\Jast^{-1}(\nabla \Jast)\circ\psi^\sharp.$$ The second identity
follows directly from \ref{sigmaom1}.

\noindent \ref{leeform}: From \ref{sigmaom1} we have
  \begin{eqnarray*}
    \sum_{A=1}^{2n}\sigma_{e_A} e_A &=&
                                        \sum_{A,B=1}^{2n}(d\omega^+(e_A,e_A, e_B) - d\omega^+(e_A,Je_A,
                                        Je_B))\; e_B\\
                                    &=&
                                        - \sum_{B=1}^{2n}\left({(\,\textstyle\sum_{A=1}^{2n} \theta_A\wedge \Jst\theta_A)}\iprod d\omega^+\right) (Je_B) \;\; e_B\\
                                    &=&
                                        - 2 \sum_{B=1}^{2n} \Jst ({\omega}\iprod d\omega^+) (e_B) \;\, e_B\ \ \ \mbox{(since
                                        $\sum_{A=1}^{2n} \theta_A\wedge \Jst\theta_A = 2\omega$)}\\
                                    &=& 
                                        -2(\Jst\theta)^\sharp.
  \end{eqnarray*}
\end{proof}

The commutator between ${\mathcal H}$ and $D_\sigma - L_{D\omega}$
will add information about commutators in the exterior algebra and will be used in the next section. We
find it here.

\begin{Lemma} \thlabel{Dsigmacomm} The following identities hold:
  \begin{enumerate}[label=\rm{(\alph*)},leftmargin=25pt]
  \item $L_{D\omega}^c = L_{D^c\omega}$.\label{Lc}
  \item $[L_{D\omega},\Hc] = i L_{\Jd D\omega}$ and\
    $[L_{D^c\omega},\Hc] = i L_{\Jd D^c\omega}$. \label{Lcomm}
  \item $D_\sigma\omega = -\Jd D\omega + 3 D^c\omega$\ and\
    $D^c_\sigma\omega = -\Jd D^c\omega - 3
    D\omega$.\label{Dsigmaomega}
  \item $[D_\sigma,\Jd] = D_\sigma^c$.\label{DsigmaJ}
  \item $[D_\sigma^c,\Jd] = -D_\sigma$.\label{DsigmacJ}
  \item
    $[D_\sigma,\Hc] = i(3 D_\sigma^c -
    L_{D_\sigma\omega})$.\label{DsigmaH}
  \item
    $[D_\sigma^c,\Hc] = -i(3 D_\sigma +
    L_{D^c_\sigma\omega})$.\label{DsigmacH}
  \item
    $[D_\sigma-L_{D\omega},\Hc] = 3i(D_\sigma^c -
    L_{D^c\omega})$.\label{DsigmaLH}
  \item
    $[D^c_\sigma-L_{D^c\omega},\Hc] = -3i(D_\sigma -
    L_{D\omega})$.\label{DsigmaLcH}
  \item $(D_\sigma-L_{D\omega})$ and $(D^c_\sigma-L_{D^c\omega})$ are
    self-adjoint.\label{DsigmaLsad}
  \end{enumerate}
\end{Lemma}

\begin{proof}
  Let $X\in \ClM$. For \ref{Lc}, note that since $D\omega$ has odd
  degree, $\Ja^{-1} D\omega = - \Ja D\omega$ by \thref{elemprop}
  \ref{2g}, and then
  $L_{D\omega}^c X = \Ja^{-1}(D\omega\cdot \Ja X) = - L_{\Ja D\omega}
  X$.  For \ref{Lcomm}, note that if $Y\in \ClM$, using the definition
  of $\Hc$ and \thref{elemprop2} \ref{wJ},
$$
[L_{Y},\Hc] X = -\dfrac{i}{2}\left(Y\cdot(\omega\cdot X +
  X\cdot\omega) - (\omega\cdot Y\cdot X + Y\cdot X\cdot\omega)\right) =
-\dfrac{i}{2}(Y\cdot\omega - \omega\cdot Y)\cdot X =iL_{\Jd Y} X.$$ We
compute $D_\sigma\omega$ directly in order to prove
\ref{Dsigmaomega}. Since $\Ja\omega=\omega$ and
$\Jd\omega = 0$, and
$\sum_{A=1}^{2n} e_A\cdot\nabla_{Je_A} X = -\sum_{A=1}^{2n}
Je_A\cdot\nabla_{e_A} X$ for any $X\in \ClM$,
\begin{eqnarray*}
  D_\sigma\omega &=& \sum_{A=1}^{2n} e_A\cdot\sigma_A(\omega)\\
                 &=& \sum_{A=1}^{2n} (e_A\cdot\nabla_{e_A}\Jd\omega-
                     e_A\cdot\Jd\nabla_{e_A}\omega + 
                     e_A\cdot\Ja^{-1}(\nabla_{Je_A} \Ja\omega) - 
                     e_A\cdot \nabla_{Je_A}\omega)  \\
                 &=& \sum_{A=1}^{2n} 
                     (- \Jd(e_A\cdot\nabla_{e_A}\omega) 
                     + Je_A\cdot\nabla_{e_A}\omega+ 
                     \Ja^{-1}(Je_A\cdot\nabla_{Je_A} \Ja\omega) +
                     Je_A\cdot\nabla_{e_A}\omega)\\
                 &=& -\Jd D\omega + D^c\omega + 2\sum_{A=1}^{2n} 
                     Je_A\cdot\nabla_{e_A}\omega \\
                 &=& -\Jd D\omega + D^c\omega + 2\sum_{A=1}^{2n} 
                     \Ja^{-1}(e_A\cdot\nabla_{e_A}\omega) \ \ 
                     \mbox{($\Ja\nabla_{e_A}\omega = -\nabla_{e_A}\omega$ by \thref{elemprop2} \ref{2n})} \\
                 &=&
                     - \Jd D\omega +3 D^c\omega 
\end{eqnarray*}
The statement about $D_\sigma^c\omega$ follows easily conjugating the
last expression with $\Ja$ and using \thref{elemprop} \ref{2h}.
   
To prove \ref{DsigmaJ}, compute
\begin{eqnarray*} [D_\sigma,\Jd] X
  &=& \sum_{A=1}^{2n} (e_A\cdot\sigma_{e_A}(\Jd X) - \Jd(e_A\cdot\sigma_{e_A} X)) \\
  &=&\sum_{A=1}^{2n} (e_A\cdot [\sigma_{e_A},\Jd] X - \Jd e_A\cdot\sigma_{e_A} X) \\
  &=&\sum_{A=1}^{2n} (-2e_A\cdot \sigma_{Je_A} X - \Jd
      e_A\cdot\sigma_{e_A} X) \ \ \mbox{by \thref{sigmaprops} \ref{sigmabra}}\\
  &=&\sum_{A=1}^{2n} Je_A\cdot \sigma_{e_A} X.
\end{eqnarray*}
On the other hand, since $\sigma_X^c = -\sigma_X$ (see
\thref{sigmaprops} \ref{sigmac}),
\begin{equation}\label{Dsigc2}
  D_\sigma^c X = \Ja^{-1}\left(\sum_{A=1}^{2n} e_A\cdot \sigma_{e_A} \Ja X\right)
  = - \sum_{A=1}^{2n} Je_A\cdot \Ja^{-1}\sigma_{e_A} \Ja X = -
  \sum_{A=1}^{2n} Je_A\cdot \sigma^c_{e_A} X =  \sum_{A=1}^{2n}
  Je_A\cdot \sigma_{e_A} X,
\end{equation}
which proves \ref{DsigmaJ}.

\ref{DsigmacJ} follows from \ref{DsigmaJ} and \thref{elemprop} \ref{2h}.

For \ref{DsigmaH}, recall (\thref{elemprop} \ref{wJ}) that
$\Hc = i(\Jd-L_\omega)$, so we only have to find
$[D_\sigma,L_{\omega}]$:
\begin{eqnarray*} [D_\sigma,L_{\omega}] X &=& \sum_{A=1}^{2n}
  (e_A\cdot
                                              \sigma_{e_A} (\omega\cdot X) - \omega\cdot \sigma_{e_A} X)\\
                                          &=& \sum_{A=1}^{2n}
                                              (e_A\cdot \sigma_{e_A}
                                              \omega\cdot X + e_A\cdot
                                              \omega\cdot\sigma_{e_A}  X- \omega\cdot e_A \cdot \sigma_{e_A} X)\\
                                          &=& L_{D_\sigma\omega} X -2
                                              \sum_{A=1}^{2n}
                                              Je_A \cdot \sigma_{e_A} X \ \ \mbox{by \thref{elemprop2} \ref{wJ}}\\
                                          &=& L_{D_\sigma\omega} X -2
                                              D_\sigma^c X \ \
                                              \mbox{by
                                              (\ref{Dsigc2}).}
\end{eqnarray*} 
Thus, $[D_\sigma, \Hc] = i(3D_\sigma^c - L_{D_\sigma\omega}).$

\ref{DsigmacH} follows from the previous one and \thref{elemprop} \ref{2h}.

\ref{DsigmaLH} and \ref{DsigmaLcH} follow immediately from
\ref{Lcomm}, \ref{Dsigmaomega}, \ref{DsigmaH}, and \ref{DsigmacH}.

To prove \ref{DsigmaLsad} we will use the fact that if $Z\in T_\C M$,
then $(L_Z)^* = -L_Z$, and if $Z\in \ClM$ is such that
$Z^\flat\in A^3_\C (M)$, then $(L_Z)^* = L_Z$ (see Propositions 9.27
and 9.29 of \cite{spinorsandcal}). If we let $X,Y\in \ClM$,
\begin{eqnarray*}
  \langle D_\sigma X,Y\rangle &=& 
                                  \sum_{A=1}^{2n}\langle e_A\cdot \sigma_{e_A} X, Y\rangle\\
                              &=& 
                                  -\sum_{A=1}^{2n}\langle  \sigma_{e_A} X, e_A\cdot Y\rangle\\
                              &=&
                                  \sum_{A=1}^{2n}\langle   X, \sigma_{e_A} (e_A\cdot Y)\rangle \ \ \mbox{by
                                  \thref{sigmaprops} \ref{sigmasa}}\\
                              &=&\sum_{A=1}^{2n}
                                  \langle   X, (\sigma_{e_A} e_A)\cdot Y + e_A\cdot \sigma_{e_A}
                                  Y\rangle 
                                  \ \ \mbox{by \thref{sigmaprops} \ref{sigmader}}\\
                              &=&
                                  \left\langle   X, \left(\sum_{A=1}^{2n}\sigma_{e_A} e_A\right)\cdot Y + D_{\sigma}
                                  Y\right\rangle\\
                              &=& 
                                  \left\langle   X, (D_\sigma - 2 L_{(\Jst\theta)^\sharp})
                                  Y\right\rangle \ \ \mbox{by \thref{sigmaprops} \ref{leeform}}.
\end{eqnarray*}
Thus,
$$D_\sigma^* = D_\sigma - 2 L_{(\Jst\theta)^\sharp}.$$
On the other hand, since $(D\omega)^\flat = d\omega + d^*\omega$ (see \thref{equiv1}\ref{Dequiv}), we can decompose $D\omega$
as $(d\omega)^\sharp+(d^*\omega)^\sharp $, and therefore
$L_{D\omega} = L_{(d\omega)^\sharp} + L_{(d^*\omega)^\sharp}$ and
$(L_{D\omega})^* = L_{(d\omega)^\sharp} - L_{(d^*\omega)^\sharp}$
since $d\omega\in A^3_\C (M)$ and $d^*\omega\in A^1_\C (M)$. Thus, using $ d^* \omega =\Jst\theta$,
$$(L_{D\omega})^* = L_{D\omega} - 2 L_{(\Jst\theta)^\sharp},$$ 
and subtracting the last displayed equations we obtain
$$(D_\sigma-L_{D\omega})^* = D_\sigma -L_{D\omega},$$
as desired. The other expression in \ref{DsigmaLsad} follows
immediately from this one and \thref{elemprop} \ref{2h}.

\end{proof}

Next, we convert the identities obtained in this section to familiar objects in the exterior algebra. Most of the proofs will involve straightforward linear algebra computations. 

\section{Conversion to the Exterior Algebra}

The formula $[D, \Hc] = -i D^c +i D_\sigma - i L_{D\om}$ proved in
Section \ref{notation} summarizes all the K\"ahler identities once it
is translated to the exterior algebra. In this section we make this
translation.

Recall that the operator $d$ decomposes into 4 parts, which we will
denote by $\mu$, $\del$, $\delb$ and $\mub$, of degrees $(2,-1)$,
$(1,0)$, $(0,1)$ and $(-1,2)$, respectively. That is,
$$d=\mu + \del + \delb + \mub.$$

From \thref{equiv1} we know that $D^c \cong -(d^c+{d^*}^c)$. 
In order to translate $D_\sigma$ and $L_{D\omega}$ we need to extend the definition of some well-known operators  in the exterior algebra (see  for example \cite{demaillytau}) to the almost-complex case. 
\begin{Definition}\thlabel{deflambda} For $\phi\in A_\C(M)$, let
  \begin{itemize}[leftmargin=20pt]\setlength\itemsep{4pt}
  \item $L\phi := \omega\wedge\phi$ (the usual Lefschetz operator),
    \  $\Lambda := L^* = \omega\iprod\phi$, \ 
    $H:=[L,\Lambda]$.
  \item
    $\lambda_\del (\phi) := \del\omega\wedge\phi$, \ 
    $\lambda_\delb (\phi) := \delb\omega\wedge\phi$ (usually denoted $\lambda$ and $\bar{\lambda}$, respectively, in the literature).\\
    $\lambda_\mu (\phi) := \mu\omega\wedge\phi$, \ 
    $\lambda_\mub (\phi) := \mub\omega\wedge\phi$.  \\
    $\lambda_+ := \lambda_\del + \lambda_\delb$,\  $\lambda_- := \lambda_\mu + \lambda_\mub$. \ \  
    $\lambda := \lambda_+ + \lambda_-$. 
  \item $\tau_\del := [\Lambda,\lambda_\del]$, \ 
    $\tau_\delb := [\Lambda,\lambda_\delb]$ (usually denoted $\tau$
    and $\overline{\tau}$, respectively, in the literature).\\
    $\tau_\mu := [\Lambda,\lambda_\mu]$, \ 
    $\tau_\mub := [\Lambda,\lambda_\mub]$.\\
   $\tau_+ := \tau_{\del} + \tau_\delb $, \
    $\tau_- := \tau_{\mu} + \tau_\mub $. \ \ $\tau = \tau_+ + \tau_-$.
    \ Note that $\tau(1) = \theta$ (the Lee form).
  \end{itemize}
\end{Definition}

To translate $L_{D\omega}$ to an operator in the exterior bundle, first note that $D\omega\cong d\omega + d^c\omega\in A^1_\C(M)\oplus A^3_\C(M)$. If $\alpha\in A^1(M)$, it is well known that
\begin{equation}
  (\alpha^\sharp\cdot\phi^\sharp)^\flat = \alpha\wedge \phi - \alpha^\sharp\iprod \phi.\label{dottowed}
\end{equation}
The equivalent statement for $\xi\in A^3(M)$ is as follows.
\begin{Lemma}\thlabel{clifmult}
  Let $\xi\in A^3_\C(M)$. For $\phi\in A_{\C}(M)$, let
    $\DS r_{\xi} (\phi):=\sum_{A=1}^{2n} (e_A\iprod \xi^\flat)\wedge
    (e_A\iprod \phi).$
Then
  \begin{enumerate}[label=\rm{(\alph*)},leftmargin=25pt]
  \item 
$(\xi^\sharp\cdot\phi^\sharp)^\flat = \xi\wedge\phi + r_{\xi}(\phi) +
(r_{\xi})^*(\phi) + \xi\iprod\phi.$ \label{clifmult1}
\item If ${\xi}$ has bidegree $(r,s)$, then the operator $r_{\xi}$ has
  bidegree $(r-1, s-1)$.  \label{clifmult2}
\end{enumerate}
\end{Lemma}
\begin{proof}
  To prove part \ref{clifmult1} we can assume, by linearity, that
  $\phi\in A^k_\C(M)$ for some $k$ and that
  $\xi=\theta^A\wedge \theta^B\wedge \theta^C$, with $A<B<C$ (so
  $\xi^\sharp=e_A\cdot e_B\cdot e_C$).  Using the formula (\ref{dottowed}) repeatedly
  it is easy to see that the operator
  $\flat\circ L_{\xi^\sharp}\circ\sharp$ has 4 components, of degrees
  $-3, -1, 1,$ and $3$ respectively.  On the other hand,
  $\flat\circ L_{\xi^\sharp}\circ\sharp$ is self-adjoint, which
  follows immediately from Propositions 9.27 and 9.29 of
  \cite{spinorsandcal}. This implies that the components of degree
  $-3$ and $-1$ of $\flat\circ L_\xi\circ\sharp$ are the adjoints of
  the components of degree $3$ and $1$, respectively.

  Thus, we only need to compute the components of degrees 3 and
  1. Using the formula (\ref{dottowed}) three times, one obtains
  \begin{eqnarray*}
    (e_A\cdot e_B\cdot e_C\cdot \phi^\sharp)^\flat &=&
                                                       \theta^A\wedge\theta^B\wedge\theta^C\wedge\phi
                                                       +\theta^A\wed\theta^C\wed (e_B\iprod \phi) - 
                                                       \theta^A\wed\theta^B\wed (e_C\iprod \phi)
                                                       -\theta^B\wed\theta^C\wed (e_A\iprod \phi)\\
                                                   &&\qquad + \ \mbox{the adjoint of the previous terms}\\
                                                   &=& 
                                                       \xi\wedge\phi
                                                       -(e_B\iprod \xi)\wed (e_B\iprod \phi) - 
                                                       (e_C\iprod \xi)\wed (e_C\iprod \phi)
                                                       -(e_A\iprod \xi)\wed (e_A\iprod \phi)\\
                                                   &&\qquad + \ \mbox{the adjoint of the previous terms}\\
                                                   &= & \xi\wedge\phi + r_{\xi}(\phi) + r_{\xi}^*(\phi) + \xi\iprod\phi.
  \end{eqnarray*}

  To prove part \ref{clifmult2}, for $1\le j\le n$, let
$$\ep_j := \dfrac{e_j - i e_{j+n}}{2} \ \ \mbox{and} \ \ \bar{\ep}_j
:= \dfrac{e_j + i e_{j+n}}{2} \ \ \mbox{so that} \ \ e_j =
\ep_j+\bar{\ep}_j \ \ \mbox{and} \ \ e_{j+n} = i(\ep_j-\bar{\ep}_j).$$
Then $\ep_j$ and $\bar{\ep}_j$ are $(1,0)$ and $(0,1)$ vectors,
respectively, and if $\psi\in A^{p,q}_\C(M)$,
\begin{eqnarray*}
  r_{\xi} (\psi) &=& \sum_{A=1}^{2n}(e_A\iprod \xi)\wed
                     (e_A\iprod \psi)\\
                 &=& \sum_{j=1}^{n}\biggl( ((\ep_j+\bar{\ep}_j)\iprod \xi)\wed
                     ((\ep_j+\bar{\ep}_j)\iprod \psi) - 
                     ((\ep_j-\bar{\ep}_j)\iprod \xi)\wed
                     ((\ep_j-\bar{\ep}_j)\iprod \psi)
                     \biggr)\\
                 &=&2 \sum_{j=1}^{n}\biggl( 
                     (\ep_j\iprod \xi)\wed (\bar{\ep}_j\iprod \psi) +
                     (\bar{\ep}_j\iprod \xi)\wed ({\ep}_j\iprod \psi)\biggr),
\end{eqnarray*}
which has bidegree $(p+r-1, q+s-1)$. Therefore, $r_{\xi}$ has bidegree
$(r-1,s-1)$.

\end{proof}
In view of this lemma, the definition of the following operator $\rho$ becomes quite natural.
\begin{Definition} For $\phi\in A_\C(M)$, let\vspace{-10pt}
  \begin{eqnarray*}
    \rho_\del(\phi):=- \sum_{A=1}^{2n} (e_A\iprod \del\omega)\wedge (e_A\iprod \phi),
    & &  \rho_\delb(\phi):=-
      \sum_{A=1}^{2n} (e_A\iprod \delb\omega)\wedge (e_A\iprod \phi),\\
  \rho_\mu(\phi):= -
\sum_{A=1}^{2n} (e_A\iprod \mu\omega)\wedge (e_A\iprod \phi), &&
 \rho_\mub(\phi):=- 
\sum_{A=1}^{2n} (e_A\iprod \mub\omega)\wedge (e_A\iprod \phi).\\ \noalign{\vspace{4pt}}
\qquad \rho_+ := \rho_\del + \rho_\delb,\qquad
\rho_- := \rho_\mu + \rho_\mub; &&
 \rho := \rho_+ + \rho_-.
  \end{eqnarray*}
\end{Definition}

With these definitions, the following proposition is almost trivial. For convenience, we use $E_\phi$ to denote exterior multiplication by $\phi$ and $I_\phi$ to denote interior multiplication by $\phi^\sharp$ (i.e. the adjoint of $E_\phi$).
\begin{Proposition}\thlabel{Domega}
  $L_{D\om}\cong \lambda + \rho + E_{\Jst\theta} +\rho^* + \lambda^* -
  I_{\Jst\theta}.$
\end{Proposition}
\begin{proof}
  Since $(D\om)^\flat = d\om + d^*\om$, with
  $d\om\in A^3_{\C}(M)$ and $d^*\om\in A^1_{\C}(M)$,
  \thref{clifmult} and formula (\ref{dottowed}) give, for
  $\phi\in A_\C (M)$,
  \begin{eqnarray*}
    \flat\circ L_{D\om}\circ\sharp(\phi) &=&(D\om\cdot\phi^\sharp)^\flat\\
                                         &=& ((d\om)^\sharp \cdot\phi + (d^*\om)^\sharp
                                             \cdot\phi)^\flat\\
                                         &=& d\om\wed\phi + r_{d\om}(\phi) + (r_{d\om})^*(\phi) +d\om\iprod\phi
                                             + d^*\om\wedge \phi - d^*\om\iprod\phi\\
                                         &=& \lambda\phi +\rho\phi + \lambda^*\phi +\rho^*\phi
                                             +E_{\Jst\theta}\phi - I_{\Jst\theta}\phi.
  \end{eqnarray*}
\end{proof}

The operators $\rho$, while quite natural, seem a bit misterious. At first it seems that they can be written in terms of the other operators in \thref{deflambda}, and in fact this is partly the case: it turns out that $\tau_\mu = i\rho_\mu$. To show this,
we first prove a general statement that will also be used later.
\begin{Lemma}\thlabel{Lambdaxi} Let $\xi$ be an odd-degree form and
  $\psi$ any form. Then
$$[\Lambda,L_\xi](\psi) = L_{\Lambda\xi}\psi + \sum_{C=1}^{2n} (e_C\iprod
\xi)\wedge (Je_C\iprod \psi).$$ In particular,
$[\Lambda,L_\xi] - L_{\Lambda\xi}$ is a derivation.
\end{Lemma}
\begin{proof} Using an adapted basis $\{e_A\}_{A=1}^{2n}$, 
  \begin{eqnarray*}
    \lefteqn{[\Lambda, L_\xi](\psi)
    = 
    \omega\iprod (\xi\wed\psi) - \xi\wed (\omega\iprod\psi)}\\ 
    &=&
        \sum_{i=1}^n \bigl( e_{i+n}\iprod  e_i\iprod (\xi\wed\psi)\bigr) - \xi\wed (\omega\iprod\psi)\\ 
    &=&
        \sum_{i=1}^n\bigl( e_{i+n}\iprod ((e_i\iprod \xi)\wed\psi)  -
        e_{i+n}\iprod  (\xi\wed (e_i\iprod\psi)) - \xi\wed (e_{i+1}\iprod
        e_i\iprod\psi)\bigr)\\ 
    &=&
        \sum_{i=1}^n\bigl( e_{i+n}\iprod e_i\iprod \xi)\wed\psi)  + 
        (e_i\iprod \xi)\wed(e_{i+n}\iprod \psi) -
        (e_{i+n}\iprod \xi)\wed (e_i\iprod\psi)\bigr)\\
    &=&  
        \Lambda\xi\wedge \psi + \DS\sum_{C=1}^{2n} (e_C\iprod \xi)
        \wedge (Je_C\iprod\psi),
  \end{eqnarray*}
  Then $[\Lambda,L_\xi] - L_{\Lambda\xi}$ is a derivation since
  $\iprod$ is and $e_C\iprod \xi$ has even degree.

\end{proof}

\begin{Lemma}\thlabel{rhotau} $\tau_\mu = i\rho_\mu$ (and
  $\tau_\mub = -i\rho_\mub$). Therefore, $\rho_{-} = -\tau_-^c$.
\end{Lemma}
\begin{proof}
  Using $\xi=\mu\omega$ in \thref{Lambdaxi}, and since
  $\Lambda(\mu\omega)=0$ because $\mu\omega$ has bidegree $(3,0)$ and
  $\Lambda$ has bidegree $(-1,-1)$, we find that
$$\tau_\mu(\psi) = \DS\sum_{C=1}^{2n} (e_C\iprod \mu\omega)
\wedge (Je_C\iprod\psi).$$ Hence $\tau_\mu$ is a derivation of degree
1, and therefore it suffices to prove the identity on scalars and
1-forms. For scalars, it is clear that both sides are 0. If $\al$ is a
1-form,
$$\Jast(\tau_\mu(\al)) =\Jast\left( \DS\sum_{C=1}^{2n} (e_C\iprod \mu\omega)\;
  \al(Je_C)\right) = - \DS\sum_{C=1}^{2n} (e_C\iprod \mu\omega)\;
\Jst\al(e_C) = \rho_\mu(\Jst\al),$$ where we have used that
$\Jast(e_C\iprod\mu\omega) = -e_C\iprod\mu\omega$ because $\mu\omega$ is
a $(3,0)$-form and therefore $e_C\iprod\mu\omega$ is a
$(2,0)$-form. Thus, $\tau_\mu = \rho_\mu^c = i\rho_\mu$ by \thref{elemprop} \ref{Pc} since $\rho_\mu$ has bidegree $(2,-1)$ (see
\thref{clifmult} \ref{clifmult2}). This also gives $\tamub = -i\rhmub$ and therefore 
 $\tau_-^c = i\tau_\mu - i \tau_\mub = -\rho_- .$

\end{proof}

In order to translate the operator $D_\sigma$ into the exterior
algebra in terms of $\rho$ and $\tau$ we need the following technical
lemmas.

\begin{Lemma} \thlabel{thtaucal}The operators $\tau_+$ and $\tau_+^c$
  can be expressed locally as follows:
  \begin{enumerate}[label=\rm{(\alph*)},leftmargin=25pt]
  \item If $\psi$ is any form,
    $\tau_+(\psi) = \theta\wedge \psi + \DS\sum_{C=1}^{2n} (e_C\iprod
    d\omega^+) \wedge (Je_C\iprod\psi)$, where $\theta$ is the Lee
    form.\label{tau+}
  \item If $\al$ is a 1-form,
    $\tau_+^c(\al) = -J\theta\wedge \al +
    \dfrac{1}{2}\DS\sum_{A,B,C=1}^{2n} d\omega^+(Je_A,e_C, Je_B)\,
    \al(e_C)\;\theta^A\wed\theta^B$.\label{taucal}
  \end{enumerate}
\end{Lemma}
\begin{proof} For part \ref{tau+}, note that $\tau_+ = [\Lambda,\lambda_+] = [\Lambda,L_{d\omega^+}]$ and use \thref{Lambdaxi}, observing that $\Lambda {d\omega^+} = \theta$ by definition.
  
In particular, when $\al$ is a 1-form,
$$\tau_+(\al) = \theta\wedge \al + \dfrac{1}{2}\DS\sum_{A,B,C=1}^{2n}
d\omega^+(e_C, e_A, e_B)\, \al(Je_C)\;\theta^A\wed\theta^B.$$
Using this, part \ref{taucal} is then a simple computation using that $\tau_+^c := \Ja^{-1}\circ\tau_+\circ \Ja$.

\end{proof}

\begin{Lemma}\thlabel{equiv2} Let
$\DS
 D_{\sigma^\flat}^{\rm ext}\psi := \sum_{A=1}^{2n} \theta^A\wedge \sigma^\flat_{e_A} \psi$
    and
    $D_{\sigma^\flat}^{\rm int}\psi := \sum_{A=1}^{2n} e_A\iprod \sigma^\flat_{e_A} \psi $.
Then
  \begin{itemize}
  \item
    $D_{\sigma^\flat}^{\rm ext} = \rho_+ + \tau_+^c + E_{\Jst\theta}$
  \item
    $D_\sigma^{\rm int} = -(D_\sigma^{\rm ext})^* + 2I_{\Jst\theta}.$
  \end{itemize}
\end{Lemma}
\begin{proof}
  For the first statement, since $\sigma^\flat_{e_A}$ is a derivation of degree 0 by \thref{sigmaprops} \ref{sigmaflatder},
  $D_{\sigma^\flat}^{\rm ext}$ is a derivation of degree 1.
  On the other hand, $\rho_+ + \tau_+^c + E_{\Jst\theta}$ is also a
  derivation of degree 1 since $\rho_+$ and
  $\tau_+^c + E_{\Jst\theta}$ are.

  Thus, it suffices to check the
  identity on 0- and 1-forms. Both sides vanish for 0-forms.
  If $\al$ is a 1-form let us use \thref{sigmaprops}
  \ref{sigmaflatder} to find an expression for
  $D_{\sigma^\flat}^{\rm ext}\al$. Using \thref{psi+l} \ref{psi++} we
  have
  \begin{eqnarray*}
    D_{\sigma^\flat}^{\rm ext}\al &=&
                                      \sum_{A,B,C=1}^{2n}  \left(d\omega^+(e_A,e_C,e_B)\,\al(e_C)\;\theta^A\wed\theta^B -
                                      d\omega^+(e_A,Je_C,Je_B)\,\al(e_C)\;\theta^A\wed\theta^B\right)\\
                                  &=&
                                      \dfrac{1}{2}\sum_{A,B,C=1}^{2n}  \left(d\omega^+(e_A,e_C,e_B)\,\al(e_C)\;\theta^A\wed\theta^B +
                                      d\omega^+(Je_A,e_C,Je_B)\,\al(e_C)\;\theta^A\wed\theta^B\right)
  \end{eqnarray*}
  The first term can be written as
$$\sum_{A,B,C=1}^{2n}
d\omega^+(e_A,e_C,e_B)\,\al(e_C)\;\theta^A\wed\theta^B =
-2\sum_{C=1}^{2n} (e_C\iprod d\omega^+)\wed (e_C\iprod \al) =
2\rho_+(\al).$$ From \thref{thtaucal} \ref{taucal}, the second term is
$$\sum_{A,B,C=1}^{2n}d\omega^+(Je_A,e_C,Je_B)\,\al(e_C)\;\theta^A\wed\theta^B
= 2(\tau_+^c (\al) + E_{\Jst\theta}(\al)).$$ Thus,
$$D_{\sigma^\flat}^{\rm ext}\al = \rho_+(\al) + \tau_+^c (\al) +
E_{\Jst\theta}(\al),$$ as desired.

For the second statement, we calculate
the adjoint directly. Let
$\phi, \psi\in A_\C^k(M)$. Then
\begin{eqnarray*}
  \left\langle D_{\sigma^\flat}^{\rm int} \phi, \psi\right\rangle &=& 
                                                                      \left\langle \sum_{A=1}^{2n} e_A\iprod \sigma^\flat_{e_A}(\phi), \psi\right\rangle\\
                                                                  &=&
                                                                      \sum_{A=1}^{2n} \langle  \sigma^\flat_{e_A}(\phi), \theta^A\wedge \psi\rangle\\
                                                                  &=&
                                                                      -\sum_{A=1}^{2n}  \langle  \phi, \sigma^\flat_{e_A} (\theta^A\wedge
                                                                      \psi)\rangle\ \ \mbox{by \thref{sigmaprops} \ref{sigmasa}}\\
                                                                  &=&
                                                                      - \sum_{A=1}^{2n}\langle   \phi, \sigma^\flat_{e_A} (\theta^A)\wedge
                                                                      \psi + \theta^A\wedge
                                                                      \sigma^\flat_{e_A}(\psi)\rangle\ \ \mbox{by \thref{sigmaprops2} \ref{sigmaflatder}}\\
                                                                  &=&
                                                                      \left\langle   \phi, 2 \Jst\theta\wedge
                                                                      \psi - D_{\sigma^\flat}^{\rm ext}(\psi)\right\rangle \ \
                                                                      \mbox{by \thref{sigmaprops2} \ref{leeform}.}
\end{eqnarray*}
Therefore,
$D_\sigma^{\rm int} = (-D_{\sigma^\flat}^{\rm ext} + 2
E_{\Jst\theta})^* = - (D_{\sigma^\flat}^{\rm ext})^* + 2
I_{\Jst\theta}$.

\end{proof}

\begin{Proposition}\thlabel{Dsigma}
  $D_\sigma\cong \rho_+ +\tau_+^c + E_{\Jst\theta} +\rho_+^* +
  {\tau_+^c}^* - I_{\Jst\theta},$ and therefore
$$D_\sigma - L_{D\om}\cong \tau^c  -\lambda  + {\tau^c}^*  -\lambda^* \ \ \ \mbox{and} \ \ \ 
D^c_\sigma - L^c_{D\om}\cong \tau  +\lambda^c  + {\tau}^*  +{
 \lambda^*}^c.$$ 
\end{Proposition}
\begin{proof}
  Using the definitions and the expression of the Clifford
  multiplication in terms of exterior and interior product, we have
$$D_\sigma \psi^\sharp=\sum_{A=1}^{2n} e_A\cdot \sigma_{e_A}\psi^\sharp \cong 
\sum_{A=1}^{2n} \theta^A\wed \sigma^\flat_{e_A}\psi -\sum_A
\theta^A\iprod \sigma^\flat_{e_A}\psi = D_{\sigma^\flat}^{\rm ext}\psi
- D_{\sigma^\flat}^{\rm int}\psi .$$ Using \thref{equiv2} we find
\begin{eqnarray*}
  D_\sigma \psi^\sharp &\cong & D_{\sigma^\flat}^{\rm ext} +
                                (D_{\sigma^\flat}^{\rm ext})^* - 2 I_{\Jst\theta}\\
                       &=& \rho_+ + \tau_+^c + E_{\Jst\theta} + \rho_+^* + {\tau_+^c}^* +
                           I_{\Jst\theta}- 2 I_{\Jst\theta}\\
                       &=&
                           \rho_+ + \tau_+^c + E_{\Jst\theta} + \rho_+^* + {\tau_+^c}^* -
                           I_{\Jst\theta}
\end{eqnarray*}
as claimed.

Since   $L_{D\om}\cong \lambda + \rho + E_{\Jst\theta} +\rho^* + \lambda^* -
  I_{\Jst\theta}$
  (\thref{Domega}), and since $\rho = \rho_+ -\tau^c_-$ by \thref{rhotau}, we obtain the second expression. For the last expression, use \thref{elemprop} \ref{2h} and \ref{Pequiv}.

\end{proof}

Putting together all the construction in this section, together with \thref{Dcomm}, we immediately obtain the result we want.

\begin{Theorem}\thlabel{main1} {\rm (K\"ahler identities for almost complex manifolds)}
  $$[d+d^*, \Lambda - L] = d^c + \tau^c -\lambda + {d^c}^* 
  + {\tau^c}^* -\lambda^* .$$
  Or, expanding and collecting
  degrees,
$$  \begin{array}{lll}
      [d, L] = \lambda & & [d, \Lambda] = {d^c}^*+{\tau^c}^*  \\ \noalign{\vspace{2pt}}  
{[d^*,L ]} = -{d^c} - {\tau^c}  &&  [d^*,\Lambda] = -\lambda^*\\
\end{array}
$$

\end{Theorem}

\begin{proof}
  Immediate from \thref{Dcomm} and \thref{equiv1,Dsigma}.
\end{proof}

\begin{Remark}
For $M$ a {\em complex} manifold (i.e. integrable), note that these are exactly the identities in \cite{demaillytau}, Th\'{e}or\`{e}me 1.1.   
\end{Remark}

\begin{Theorem}\thlabel{main2} The following identity holds:
$$[\tau^c  -\lambda  +{\tau^c}^*  -\lambda^*, \Lambda - L] = 
3( \tau +\lambda^c +{\tau}^*
+{\lambda^c}^*).$$ Or, expanding and collecting degrees,
$$  \begin{array}{lll}
[\lambda, L] = 0 & \qquad& [\lambda, \Lambda] = -\tau \\ \noalign{\vspace{2pt}}  
[\tau,L] = -3\lambda  & \qquad& [\tau,\Lambda] = -2 {\tau^c}^*
\end{array}
$$

\end{Theorem}
\begin{proof}
  Follows immediately from \thref{equiv1} \ref{Hcequiv},
  Lemma \thref{DsigmaLH} and \thref{Dsigma}.

\end{proof}

\begin{Corollary} The following identities, together with their adjoints and
  conjugates, hold:
  \begin{equation}\label{genKI}
    \begin{array}{lll}
      [\mu,\Lambda] = i(\mub^* + \tau_\mub^*)
      &\qquad[\tau_\mu,\Lambda] =-2i  \tau_\mub^*
      &\qquad \left[\lambda_\mu,\Lambda\right] = -\tau_\mu\\
      {[\mu,L]} = \lambda_{\mu}
      &\qquad [\tau_\mu,L] = -3 \lambda_\mu
      &\qquad \left[\lambda_\mu,L\right] = 0\\
      \left[\del, \Lambda\right] = -i({\delb^*} +\tau_{\delb}^*)  
      &\qquad[\tau_\del, \Lambda ] = 2i\tau_{\delb}^*
      &\qquad[\lambda_\del,\Lambda] = -\tau_\del\\
        {[\del,L]} = \lambda_\del
      &\qquad [\tau_\del,L] = -3\lambda_\del
       &\qquad \left[\lambda_\del,L\right] = 0
    \end{array}
  \end{equation}

\end{Corollary}

\begin{proof}
All the identities follow from \thref{main1} and \thref{main2} by separating by bidegrees.

\end{proof}

\begin{Remark}
  If $M$ is almost K\"ahler (that is, $d\omega=0$), then
  $\rho, \lambda$ and $\tau$ are all $0$, and (\ref{genKI}) reduce to
$$
\begin{array}{ll} [\mu,\Lambda] = i\mub^*
  & \qquad  [\mu,L] = 0\qquad\qquad \\
  \left[ \del , \Lambda \right] = -i{\delb^*} & \qquad [\del,L] = 0
\end{array}
$$
(together with their adjoints and conjugates). These identities appear
in \cite{ciriciwilson}, Proposition 3.1, and \cite{bartolo}, p. 1295.
\end{Remark}

For reference, we give the following table of commutators of the common operators with $L$, $\La$, and $H$.

\begin{figure}[!h]
  \begin{minipage}{0.4\textwidth}
$$
\begin{array}{!{\vline width .8pt}c!{\vline width .8pt}c|c|c!{\vline width .8pt}}
  \hline
  \noalign{\vskip-.07pt}\hline
  \grcel & \La & L \\
  \hline
  \noalign{\vskip-.07pt}\hline
  d & {d^c}^*+{\tau^c}^* & \lambda \\
  \hline
  \mu & i(\mub^*+\tamub^*) & \lamu \\
  \hline
  \tamu & -2i\tamub^* & -3\lamu \\
  \hline
  \mub & -i(\mu^*+\tamu^*) & \lamub \\
  \hline
  \tamub & 2i\tamu^* & -3\lamub \\
  \hline
  \del & -i(\delb^*+\tadeb^*) & \lade \\
  \hline
  \tade & 2 i \tadeb^* & -3\lade \\
  \hline
  \rhde & -i\rhdeb^*+\tadeb^* & i\lade \\
  \hline
  \delb & i(\del^*+\tade^*) & \ladeb \\
  \hline
  \tadeb & -2 i \tade^* & -3\ladeb \\
  \hline
  \rhdeb & i\rhde^*+\tade^* & -i\ladeb \\
  \hline
  \lamub & -\tamub & 0 \\
  \hline
  \ladeb & -\tadeb & 0 \\
  \hline
  \lade & -\tade & 0 \\
  \hline
  \lamu & -\tamu & 0 \\
  \hline
  \noalign{\vskip-.07pt}\hline
\end{array}
$$
\end{minipage}
\begin{minipage}{0.4\textwidth}
  \[
    \begin{array}{!{\vline width .8pt}c!{\vline width
          .8pt}c|c|c!{\vline width .8pt}}
      \hline
      \noalign{\vskip-.07pt}\hline
      \grcel & \La & L \\
      \hline
      \noalign{\vskip-.07pt}\hline
      d^* &-\lambda^* & -({d^c}+{\tau^c}) \\
      \hline
      \mu^* &i(\mub+\tamub) &  -\lamu^* \\
      \hline
      \tamu^* &-2i\tamub &  3\lamu^* \\
      \hline
      \mub^* &-i(\mu+\tamu) &  -\lamub^* \\
      \hline
      \tamub^* &2i\tamu &  3\lamub^* \\
      \hline
      \del^* &-i(\delb+\tadeb) &  -\lade^* \\
      \hline
      \tade^* &2 i \tadeb &  3\lade^* \\
      \hline
      \rhde^* & -i\rhdeb-\tadeb   &  i\lade^* \\
      \hline
      \delb^* &i(\del+\tade) &  -\ladeb^* \\
      \hline
      \tadeb^* &-2 i \tade &  3\ladeb^* \\
      \hline
      \rhdeb^* &i\rhde-\tade  &  -i\lade^* \\
      \hline
      \lamub^* &\tamub^* &  0 \\
      \hline
      \ladeb^* &\tadeb^* &  0 \\
      \hline
      \lade^* &\tade^* &  0 \\
      \hline
      \lamu^* &-i\tamu^* &  0\\
      \hline
      \noalign{\vskip-.07pt}\hline
    \end{array}
  \]
\end{minipage}

\caption{Full table of commutators with $L$ and $\Lambda$.}
\end{figure}

\begin{Remark}
  The commutators $[\rho_\del, L] = i\lambda_\del$ and
  $[\rho_\del, \Lambda] = -i \rho_\delb^* +\tau_\delb^*$ on the
  table above do not follow immediately from the results above, but we
  include them for reference.  They were found using the methods
  developed in \cite{diracformula, samsthesis}.

\end{Remark}

\newpage

\section{Appendix}\label{appendix}

\begin{figure}[H]
  \begin{tabular} { !{\vline width .8pt} m{.3in} !{\vline width .8pt}
      m{.5in} | m{.5in}| m{.5in}| m{.5in}| >{\columncolor{gray!30}}
      m{.5in}| m{.5in}| m{.5in}| m{.5in}| m{.5in} !{\vline width
        .8pt}} \hline
    \thead{\raise-3pt\hbox{$\boldsymbol{q}$}$\backslash$
      \raise4pt\hbox{$\!\!\boldsymbol{p}$}} &
    \thead{$\boldsymbol{-4}$} & \thead{$\boldsymbol{-3}$} &
    \thead{$\boldsymbol{-2}$} & \thead{$\boldsymbol{-1}$} &
    \cellcolor{white} \thead{$\boldsymbol{0}$} & \thead{1} & \thead{2}
    & \thead{3} & \thead{4}\\
    \hline\noalign{\vspace{-.8pt}}\hline\noalign{\vspace{.4pt}}
    \ \ $\boldsymbol{4}$ &&&&&& $[\lambda_\mub,L]$&&\phantom{$A^3_g$}&\\ \hline
    \ \ $\boldsymbol{3}$ &&&&&$\lambda_\mub$\quad $[{\mub},L]$ $[\tau_{\mub},L]$ &&
                                                                                    $[\lambda_\delb,L]$&&\\ \hline 
    \ \ $\boldsymbol{2}$ &&&&$\mub$  \  $\tau_\mub$  

                              $[\mu^*,L]$ $[\tau_\mu^*,L]$
                              $[\lambda_\mub,\Lambda]$&& $\lambda_\delb$\quad $[\delb,L]$ $[\tau_\delb,L]$ && 
                                                                                                              $[\lambda_\del,L]$&\\ \hline 
    \ \ $\boldsymbol{1}$ &&&$\mu^*$ $\tau_\mu^*$ 

                             $[\lambda_\mu^*,L]$ $[\mub,\Lambda]$ 
                             $[\tau_\mub,\Lambda]$&&$\delb\,\,\,\tau_\delb$ $[\del^*,L]$ $[\tau_\del^*,L]$ $[\lambda_\delb,\Lambda]$
                              &
                                                          &$\lambda_\del$\quad $[\del,L]$ $[\tau_\del,L]$&
                                                           &$[\lambda_\mu,L]$\\ \hline \rowcolor{gray!30}
    \ \ \cellcolor{white}$\boldsymbol{0}$ &&$\lambda_\mu^*$
                                             $[\mu^*,\Lambda]$ 

                                             $[\tau_\mu^*,\Lambda]$&&$\del^*\,\,\tau_\del^*$ $[\lambda_\del^*,L]$
                                                                      $[\delb,\Lambda]$ $[\tau_\delb,\Lambda]$&&$\del\,\,\,\tau_\del$ $[\delb^*,L]$ $[\tau_\delb^*,L]$ $[\lambda_\del,\Lambda]$&&
                                                                                                                                                                                                  $\lambda_\mu$ $[\mu,L]$ $[\tau_\mu,L]$&\\ \hline 
    $\boldsymbol{-1}$ &$[\lambda_\mu^*,\Lambda]$&&$\lambda_\del^*$
                                                   $[\del^*,\Lambda]$ $[\tau_\del^*,\Lambda]$&&$\delb^*\,\,\tau_\delb^*$
                                                                                                $[\lambda_\delb^*,L]$ $[\del,\Lambda]$ $[\tau_\del,\Lambda]$& 
                                                          &$\mu$\! $\tau_\mu$ $[\mub^*, L]$ $[\tau_\mub^*, L]$ $[\lambda_\mu, \Lambda]$&
                                                           &\\ \hline 
    $\boldsymbol{-2}$ &&$[\lambda_\del^*,\Lambda]$&&$\lambda_\delb^*$
                                                     $[\delb^*,\Lambda]$ $[\tau_\delb^*,\Lambda]$&&
                                                                                                    $\mub^*$ $\tau_\mub^*$ 

                                                                                                    $[\lambda_\mub^*,L]$ $[\mu,\Lambda]$ $[\tau_\mu,\Lambda]$&&
                                                          &\\ \hline 
    $\boldsymbol{-3}$ &&&$[\lambda_\delb^*,\Lambda]$&&$\lambda_\mub^*$
                                                       $[\mub^*,\Lambda]$ $[\tau_\mub^*,\Lambda]$& 
                                                          &&
                                               &\\ \hline 
    $\boldsymbol{-4}$ &&&&$[\lambda_\mub^*,\Lambda]$ &\phantom{$A^3_g$}&&&
                                               &\\ \noalign{\vspace{.4pt}}\hline\noalign{\vspace{-.8pt}}\hline 
  \end{tabular}
  \vspace{.2in}
  \caption{Table of bidegrees $(p,q)$ of operators}
\end{figure}

\newpage
\bibliographystyle{amsplain} \bibliography{cliffordbib}{}
\providecommand{\bysame}{\leavevmode\hbox to3em{\hrulefill}\thinspace}
\providecommand{\MR}{\relax\ifhmode\unskip\space\fi MR }
\providecommand{\MRhref}[2]{%
  \href{http://www.ams.org/mathscinet-getitem?mr=#1}{#2} }
\providecommand{\href}[2]{#2}

\vskip.3in

\hbox{
\qquad\vtop{\hsize3in\baselineskip10pt\parindent0pt
Samuel Hosmer

samuel.hosmer@gmail.com

}
\vtop{\hsize3in\baselineskip10pt\parindent0pt
Luis Fern\' andez

Department of Mathematics

CUNY Bronx CC and CUNY Graduate Center

365 5$^{\rm th}$ Avenue

New York, NY 10016, USA

{\it lmfernand@gmail.com}

{\it luis.fernandez01@bcc.cuny.edu}

{\it lfernandez1@gc.cuny.edu}

}
}

\end{document}